\documentclass[11pt]{amsart} 
\usepackage[parfill]{parskip}   
\usepackage{amsthm}
\usepackage{textcomp}
\usepackage{blkarray}
\usepackage[scaled=.85]{beramono}% a typewriter font must be defined
\usepackage[leqno]{amsmath} %use eq no on left
\usepackage{bm}% load after all math to give access to bold math
\usepackage{url}
\arraycolsep=1.4pt\def\arraystretch{1.4}

\newcommand{\bC}{{\mathbb{C}}}

\newcommand{\bG}{{\mathbb{G}}}

\usepackage[paperwidth=7in, paperheight=10in, textwidth=4.75in, textheight=8in, includehead=true, bindingoffset=.5in]{geometry}

\newcommand{\beq}{\begin{equation}}
\newcommand{\eeq}{\end{equation}}
\numberwithin{equation}{section}

\newtheorem{theorem}[equation]{Theorem}

\newtheorem{lemma}[equation]{Lemma}

\newtheorem{corollary}[equation]{Corollary}

\theoremstyle{remark}

\newtheorem{remark}[equation]{Remark}

\theoremstyle{definition}

\newtheorem{definition}[equation]{Definition}

\title{\textsf{Schur Inequalities}}%\\[15pt]{\em proof techniques and intuition}}
\author{S. Gill Williamson}
\thanks{Department of Computer Science and Engineering, 
University of California San Diego; \url{http://cse.ucsd.edu/~gill}.
}
\date{}                                           % Activate to display a given date or no date
\usepackage[pdftex,colorlinks=true, pdfstartview=FitV, linkcolor=blue, citecolor=blue, urlcolor=blue,bookmarks,bookmarksnumbered]{hyperref}
\hypersetup{
pdftitle={Matrix Theory},
pdfauthor={S. Gill Williamson}%{\textcopyright\ S. Gill Williamson}
}
\begin{document}
\thispagestyle{empty}
\begin{center}
\vspace*{1in}
\textsf{ \Large Recursive projections of symmetric tensors\\
	and Marcus's proof of the Schur inequality}\\[.2in]
%\textsf{\small basic skills and proof techniques}\\[.7in]

\textsf{S. Gill Williamson}\footnote{\url{http://cseweb.ucsd.edu/~gill}}
%\vfill
%\textcopyright  \textsf{S. Gill Williamson 2012. All rights reserved.}
\end{center}

%\newpageI
\thispagestyle{empty}
\hspace{1 pt}
%\newpage
%ABSTRACT
\begin{center}
{\Large Abstract}\\[.2in]
\end{center}
\pagestyle{plain}

In a 1918 paper  \cite{is:fgh}, Schur proved a remarkable inequality that related group representations, Hermitian forms and  determinants.  He also gave concise necessary and sufficient  conditions for equality. 
In  \cite{mm:crs}, Marcus  gave a beautiful short proof of Schur's inequality by applying the Cauchy-Schwarz inequality to symmetric tensors, but he did not discuss the case of equality. In \cite{gw:tch}, Williamson gave an inductive proof of Schur's equality conditions by contracting Marcus's symmetric tensors onto lower dimensional subspaces where they remained symmetric tensors.  Here we unify these results notationally and conceptually, replacing contraction operators with the more geometrically intuitive projection operators. 
%INTRODUCTION
\section{Introduction}
The following theorem will be the focus of this paper.
%SCHUR'S THEOREM
\begin{theorem}[\bfseries Schur's theorem for finite groups and hermitian forms]
\label{thm:stf}
Let $G$ be a subgroup of the symmetric group $S_n$ on $\underline{n}=\{1,2, \ldots, n\}$.  
Let $H=(h_{ij})$ be an $n\times n$ complex positive definite Hermitian matrix and define
$\,\bG_H$ to be the group generated by all transpositions $(i,j)$ such that $h_{ij}\neq 0$.
Let $\,M$ be a representation of $\,G$ as unitary linear operators on $\,U$,
$\dim(U)=m$, and   
let $\,M_H = \sum_{\sigma\in G} M(\sigma) \prod_{i=1}^n h_{i\sigma(i)}$.
$\,M_H$, called a generalized matrix function, is  positive definite Hermitian, and for $\,u\in U$, $\|u\|=1$,
\begin{equation}
\label{eq:sin}
\det(H)\leq (M_Hu,u)
\end{equation}
\begin{equation}
\label{eq:seq}
\det(H) = (M_Hu,u) \;{\it if\;and\;only\;if}\;\bG_H\subseteq G\,\, {\it and}\,\,
\end{equation}
\[ 
(M(\sigma)u,u)=\epsilon(\sigma) \;{\it for\;all\;} \sigma\in \bG_H\;\;
{\it where\;\;}  \epsilon(\sigma) {\it \;\;is\; the\; sign\; of\; \sigma.}\\
\]
\end{theorem}
In this section we discuss the basics.
In Section~2 we prove Marcus's generalization of Schur's inequality.
In Section~3 we treat the case of equality.
In Section~4 we discuss the combinatorial lemmas needed for Section~3.
In Section~5 we  give examples of Schur's inequality and discuss the trace version.

%END SCHUR'S THEOREM
%REM
\begin{remark}[\bfseries Comments on theorem~\ref{thm:stf}]
\label{rem:comequcon}
\hspace*{1in}\\[4pt]
$\bf(1)$ Throughout the rest of the paper, ${\bf M}_n(\bC)$ denotes the $n\times n$ matrices with entries in $\bC$, the complex numbers.  Likewise, the order (or cardinality) of the finite group $G$ will be denoted by $g$.\\
\hspace*{1in}\\[1pt]
$\bf(2)$ Referring to equality condition~\ref{eq:seq} we have
\[M_H = \sum_{\sigma\in G} M(\sigma) \prod_{i=1}^n h_{i\sigma(i)}\;\,{\rm implies\; that\;\,}
(M_Hu,u) = \sum_{\sigma\in G} (M(\sigma)u,u) \prod_{i=1}^n h_{i\sigma(i)}.\]
Assume $M(\sigma_1)$ and $M(\sigma_2)$ unitary and $u\in U$ has $\|u\|=1$.    
If $|(M(\sigma_1)u,u)|=|(M(\sigma_2)u,u)|=1$ then $u$ is a unit eigenvector of both
$M(\sigma_1)$ and $M(\sigma_2)$.  
We have $M(\sigma_j)u=\lambda_ju$, $\lambda_j =e^{ir_j}$, $i=1,2$.
Thus, $(M(\sigma_1)u,u)(M(\sigma_2)u,u) = \lambda_1\lambda_2 = 
(M(\sigma_1)M(\sigma_2)u,u) = (M(\sigma_1\sigma_2)u,u).$\\
\hspace*{1in}\\[1pt]
$\bf(3)$ A group $\bG_D$, analogous to $\bG_H$ of theorem~\ref{thm:stf}, can defined for any $D\in {\bf M}_n(\bC).$ 
Let $K\subseteq \{1,2,\ldots,n\}$ be an orbit of $\bG_D$.
Since $\bG_D$ is generated by transpositions, the restriction $\bG_D|K$ is the symmetric group $S_K$.\\
\end{remark}
%END REM
%INTRODUCE TENSORS
%REM
\begin{remark}[\bfseries Inner products on tensors]
\label{rem:innproten}
Let ${\underline n}=\{1, \ldots, n\}$.  
Denote by ${\underline n}^{\underline n}$  the set of all functions
from ${\underline n}$ to ${\underline n}$.  
The notation $\Gamma_n \equiv {\underline n}^{\underline n}$ is also common in this subject.  Note the cardinality, $| \Gamma_n| = | {\underline n}^{\underline n}| =  n^n.$
Let $e_1, \ldots, e_n$ be an orthonormal basis for the unitary space $V$.  
Using the inner product in $U$, we have an inner product on 
$U\otimes [\otimes^n V]$ which on homogeneous elements is
\[
\left(u_1\otimes x_1\otimes\cdots\otimes x_n, u_2\otimes y_1\otimes\cdots\otimes y_n\right)=
(u_1,u_2)\prod_{i=1}^n (x_i,y_i).
\]
Note that
\begin{equation}
\label{eq:onbasis}
\{e_\alpha\mid e_\alpha=e_{\alpha(1)}\otimes \cdots \otimes e_{\alpha(n)},\,\,\alpha\in{\underline n}^{\underline n}\}
\end{equation}
is an orthonormal basis for $\otimes^n V$.\\
\end{remark}
%END REM
%DEF SYMMETRY OPERATORS
\begin{definition}[\bfseries Generalized symmetry operators]
\label{def:gsop}
Let $G$ be a subgroup of the symmetric group $S_n$ on $\{1,2, \ldots, n\}$. 
Let $M$ be a representation of $G$ as unitary linear operators on a unitary space $U$, $\dim U=m$.
Define an endomorphism of $U \otimes \left[\otimes^nV\right]$ by
$T_G=\sum_{\sigma\in G} M(\sigma)\otimes P(\sigma)$ where $P(\sigma)$ is the permutation operator 
defined by 
$P(\sigma)(e_\alpha)=e_{\alpha(\sigma^{-1}(1))}\otimes \cdots \otimes e_{\alpha(\sigma^{-1}(n))}$
on the basis and extended (conjugate) linearly to $\otimes^n V$.
On homogeneous tensors 
\begin{equation}
\label{eq:psigma}
P(\sigma)(x_1\otimes \cdots \otimes x_n) = x_{\sigma^{-1}(1)}\otimes \cdots \otimes x_{\sigma^{-1}(n)}.
\end{equation}
$T_G$ will be called a {\em generalized symmetry operator} on $U\otimes (\otimes^n V)$
of {\em degree} $m=\dim(U)$ and {\em order} $n=\dim(V)$.
On homogeneous elements
\[
T_G(u\otimes x_1\otimes \cdots \otimes x_n)=
\sum_{\sigma\in G} M(\sigma)u\otimes x_{\sigma^{-1}(1)}\otimes \cdots \otimes x_{\sigma^{-1}(n)}.\\
\]
\end{definition}
%END DEF
%\hspace*{1in}\\
\begin{remark}[\bfseries Cauchy-Schwartz on symmetric tensors]
\label{rem:cssymten}
Let $e=e_1\otimes \cdots \otimes e_n$,  $x=x_1\otimes \cdots \otimes x_n$ and $u\in U$, $||u||=1$.
By the Cauchy-Schwarz inequality,
%EQUQTION
\begin{equation}
\label{eq:csifront}
|\left(T_G(u\otimes x), T_G(u\otimes e)\right)|^2\leq \|(T_G(u\otimes x)\|^2\|(T_G(u\otimes e)\|^2.
\end{equation}
%END EQUATION
Let $A\in {\bf M}_n(\bC)$, $A=(a_{ij})$, be upper triangular and nonsingular. 
Define vectors $x_i$, $i=1,\ldots, n$,  by $x_i=\sum_{j=1}^n a_{ij}e_j$.
Thus, $A=(a_{ij})=((x_i,e_j))$.
We say that $A$ is {\em row--associated} with   the vectors
$x_i$, $i=1,\ldots, n$ and the basis $e_i$, $i=1, \ldots, n$.
Marcus (\cite{mm:crs}) observed that 
\[
0< |\left(T_G(u\otimes x), T_G(u\otimes e)\right)|^2= g^2(\det(A))^2 
\]
\[
\|(T_G(u\otimes x)\|^2 = g(M_{AA^*}u,u)\;\;{\rm and}\;\;
\|T_G(u\otimes e)\|^2=g.\\
\]
\end{remark}
From definition~\ref{def:gsop} and remarks~\ref{rem:innproten}, \ref{rem:cssymten}, we have the following:
%THM
\begin{theorem}[\bfseries Cauchy-Schwartz-Schur inequalities]
\label{thm:inequalities}
Let $u\in U$, $\|u\| = 1$.
Let $H\in{\bf M}_n(\bC)$ be positive definite Hermitian.  
Write $H=AA^*$ where $A$ is upper triangular and nonsingular (e.g., using Cholesky decomposition), and 
define  $x_i$, $i=1,\ldots, n$, by $A=((x_i,e_j))$ (remark~\ref{rem:cssymten}). The following inequalities are equivalent
\begin{description}
\item[(1) C-S]
$|\left(T_G(u\otimes x), T_G(u\otimes e)\right)|^2\leq \|T_G(u\otimes x)\|^2\|T_G(u\otimes e)\|^2$
\item[(2) Schur H]
$\det(H)\leq (M_{H}u,u)$,
\item[(3)  A form]
$(\det(A))^2 \leq (M_{AA^*}u,u)$.
\end{description}
\begin{proof}
This result follows directly from remark~\ref{rem:cssymten}.\\
\end{proof}
\end{theorem}
%END THM
Referring to theorem~\ref{thm:inequalities}, we have the following:
\begin{theorem}[\bfseries Cauchy-Schwartz-Schur equalities]
\label{thm:equalities}
Let $u\in U$ and $H=AA^*$ be as in theorem~\ref{thm:inequalities}.
The following equalities are equivalent.
\begin{description}
\item[($0$) C-S Equality condition]
$T_G(u\otimes x) = k T_G(u\otimes e),\;k\neq 0$
\item[($1$) C-S Equality]  
$|\left(T_G(u\otimes x), T(u\otimes e)\right)|^2= \|T_G(u\otimes x)\|^2\|T_G(u\otimes e)\|^2$
\item[($2$) Schur H Equality]
$\det(H) = (M_{H}u,u)$.
\item[($3$) A Equality]
$(\det(A))^2 = (M_{AA^*}u,u)$
\item[($4$) Schur H Equality Condition]
$\bG_H\subseteq G, (M(\sigma)u,u)=\epsilon(\sigma), \sigma\in \bG_H$
\item[($5$) A Equality Condition]
$\bG_A\subseteq G, (M(\sigma)u,u)=\epsilon(\sigma), \sigma\in \bG_A$.
\end{description}
\begin{proof}
The equivalence of conditions {\bf($0$)}, {\bf($1$)}, {\bf($2$)} and  {\bf($3$)} follows from 
theorem~\ref{thm:inequalities} and the standard condition for equality in  the Cauchy-Schwartz inequality (i.e. the equivalence of {\bf($0$)} and {\bf($1$)}).\\[5pt]
%is equivalent to {\bf($1$)}, and  
%{\bf($1$)} is equivalent to {\bf($2$)} and {\bf($3$)} (Marcus, \cite{mm:crs}).
%Thus, {\bf($0$)}, {\bf($1$)}, {\bf($2$)} and {\bf($3$)} are all equivalent.
Condition {\bf($5$)} implies {\bf($4$)} (\cite{gw:tch}, Lemma 5.3):\\
{\bf Assume {\bf($5$)} and show {\bf($4$)}}. 
Note that $H(i,j)\equiv h_{ij} \neq 0$  implies that 
$H(i,j)=\sum_{k=1}^n A(i,k)A^*(k,j) = \sum_{k=1}^n A(i,k)\bar{A}(j,k)\neq 0$.
Thus, there is some $s$ such that $(i,s), (j,s)\in \bG_A$ which, by hypothesis, implies that $(M((i,s))u,u)=-1$ and $(M((j,s))u,u)=-1$.  
Thus, $(i,j)=(i,s)(j,s)(i,s)\in \bG_H$ has $(M((i,j))u,u)=-1$,
and item {\bf($4$)} holds  (we use remark~\ref{rem:comequcon} (2)).
Note this also shows that $G_H\subseteq G_A$.\\[5pt]
Condition {\bf ($4$)}  implies {\bf ($2$)} (\cite{gw:tch}, Lemma 5.2): \\
{\bf Assume {\bf($4$)} and show {\bf($2$)}}.  
$M_H = \sum_{\sigma\in G} M(\sigma) \prod_{i=1}^n h_{i\sigma(i)}$, thus
$(M_Hu,u) = \sum_{\sigma\in G} (M(\sigma)u,u) \prod_{i=1}^n h_{i\sigma(i)}$.
Note that $\sigma\in S_n$, $\prod_{i=1}^n h_{i\sigma(i)}\neq 0$ implies  
$h_{i\sigma(i)}\neq 0$, $1\leq i \leq n$, which implies, by definition of $\bG_H$, 
$(i,\sigma(i))\in\bG_H$, $1\leq i \leq n$.  
Recall that $\sigma$ is a product of transpositions of the form $(i,\sigma(i))$.
Each such transposition has $(M((i,\sigma(i)))u,u)=-1$ by our hypothesis. 
Thus, $\sigma\in \bG_H\subseteq G$ and $(M(\sigma)u,u)= \epsilon(\sigma)$.
We have shown that 
\[
(M_Hu,u) = \sum_{\sigma\in \bG_H} (M(\sigma)u,u) \prod_{i=1}^n h_{i\sigma(i)}=
 \sum_{\sigma\in S_n} \epsilon(\sigma) \prod_{i=1}^n h_{i\sigma(i)}=\det(H).
\]
Finally, item {\bf($0$)} implies {\bf($5$)} is proved in \cite{gw:tch} (Theorem 4.2)
and below (Theorem~\ref{thm:mainresult}).
This completes the proof of theorem~\ref{thm:equalities}.\\
\end{proof}
\end{theorem}
%REM
\begin{remark}[\bfseries Basic observations on theorem~\ref{thm:equalities}]
\label{rem:basicstuff}
In showing $(5)$ implies $(4)$ above, we have shown that 
$\bG_H\subseteq \bG_A$. 
Using the ideas in $(4)$ implies $(2)$, for any matrix $D\in {\bf M}_n(\bC)$,
\[
\sum_{\sigma\in G} M(\sigma)\prod_{i} D(i,\sigma(i))=
\sum_{\sigma\in \bG_D\cap G}M(\sigma)\prod_{i} D(i,\sigma(i)).
\]
From theorem~\ref{thm:equalities}, we will have $(4)$ implies $(5)$.
If we take $G=\bG_H$ in $(4)$ and $M(\sigma)=\epsilon(\sigma)I_m$, we have 
$\bG_A\subseteq \bG_H$ from $(5)$.  
Reversing this argument using $(5)$ implies $(4)$ we get 
$\bG_H\subseteq \bG_A$  (also derived from the argument proving $(5)$ implies
$(4)$ above).  Thus, $\bG_H=\bG_A$. A simple example shows the generating sets $\{(i,j)\mid A(i,j)\neq 0\}$ and $\{(i,j)\mid H(i,j)\neq 0\}$ need not  be the same.
Schur in~\cite{is:fgh} proved $(2)$ of theorem~\ref{thm:inequalities}
and showed the equivalence of $(2)$ and $(4)$ of theorem~\ref{thm:equalities}.
The focus of this paper is to characterize the multilinear algebraic properties of the symmetry operators $T_G$ that result in the equivalence of  $(0)$ and $(5)$ 
in theorem~\ref{thm:equalities}.  The most interesting structural properties of $T_G$ arise in proving $(0)$ implies $(5)$ (theorem~\ref{thm:mainresult}).  
The converse $(5)$ implies $(0)$  follows from theorem~\ref{thm:equalities}.
\end{remark}

%END THEOREM

%SECTION
\section{Inequalities}

We use the terminology of the previous section.
Let $U$ and $V$ be $m$ and $n$ dimensional unitary spaces with the standard inner products.\\

\begin{remark}[\bfseries Properties of $T_G$]
\label{rem:pot}
To summarize, the $M(\sigma)$ and $P(\sigma)$ are unitary operators: $(M(\sigma))^*=M(\sigma^{-1})$ (by definition of $M$) and 
$(P(\sigma))^*=P(\sigma^{-1})$ (simple computation). Thus, $M(\sigma)\otimes P(\sigma)$ is unitary since
 \[
 (M(\sigma)\otimes P(\sigma))^* = (M(\sigma))^*\otimes (P(\sigma))^*=M(\sigma^{-1})\otimes P(\sigma^{-1}).
 \]
 $M$ is defined to be a representation of $G$  as unitary operators on $U$, and $P$ is a representation of $G$ as unitary operators on $\otimes^n V$. 
$T_G$ is essentially idempotent, $T_G^2=gT_G$, and Hermitian, $T_G^*=T_G$.\\
\end{remark}

%REM
\begin{remark}[\bfseries Inner products as products of associated matrices]
\label{rem:ippm}
Consider the unitary inner product defined by
\begin{equation}
\label{eq:uipxy}
(x_1\otimes \cdots \otimes x_n, y_1\otimes \cdots \otimes y_n)= \prod_{i=1}^n (x_i,y_i).\\
\end{equation}
Let $x_i=\sum_{t=1}^n a_{it}e_t$ and $y_j=\sum_{t=1}^n b_{jt}e_t$ where $e_1, \ldots, e_n$ 
is the orthonormal basis for $V$.
In the unitary space $V$, we have $(x_i,y_j)=\sum_{t=1}^n a_{it} \bar{b}_{jt}$.
In terms of remark~\ref{rem:cssymten}, $A$ is row--associated with   the vectors
$x_i$, $i=1,\ldots, n$ and the basis $e_i$, $i=1, \ldots, n$, and 
$B$ is row--associated with   the vectors $y_i$, $i=1,\ldots, n$ and the basis 
$e_i$, $i=1, \ldots, n$.  With the basis $e_i$, $i=1, \ldots, n$ understood, we also refer
to $A$ and $B$ as the matrices associated with 
$x= x_1\otimes \cdots \otimes x_n$ and $y_1\otimes \cdots \otimes y_n$ respectively.

% $A=(a_{ij})$ is the matrix {\em associated} with $x= x_1\otimes \cdots \otimes x_n$ with
%respect to the o.n. basis $e_1, \ldots, e_n$.  
%The row, $A_{(i)}$, is  the vector of coefficients of $x_i$ with respect to $e_1, \ldots, e_n$.
%Likewise, let $B$ be the matrix associated with $y=y_1\otimes \cdots \otimes y_n$.
For $B=(b_{ij})$, the conjugate transpose $B^*(i,j)=(\bar{b}_{ji})$.
Thus, the inner product $(x_i,y_i)$ becomes
\begin{equation}
\label{eq:ipmtrx}
(x_i,y_j)=\sum_{t=1}^n a_{it} \bar{b}_{jt}=\sum_{t=1}^n A(i,t)B^*(t,j)=AB^*(i,j).
\end{equation}
If $x= x_1\otimes \cdots \otimes x_n=y_1\otimes \cdots \otimes y_n$ then $A=B$ and
$AA^*$ is the Gram matrix of the sequence $x_1, \ldots, x_n$ (i.e. the rows of $A$).\\
\end{remark}

%LEMMA
We use the terminology of theorem~\ref{thm:stf} and equation~\ref{eq:ipmtrx}.
\begin{lemma}[\bfseries Marcus's inner product form for  symmetric tensors \cite{mm:crs}]
\label{lem:ipgst}
Let $T_G(u_1\otimes x)$ and $T_G(u_2\otimes y)$ be two generalized symmetric tensors where 
$x= x_1\otimes \cdots \otimes x_n$, $y=y_1\otimes \cdots \otimes y_n$, and
$A$ and $B$ are the matrices associated with  $x$ and $y$ respectively
(remark~\ref{rem:ippm}).
Then 
\begin{equation}
\label{eq:iptst}
(T_G(u_1\otimes x), T_G(u_2\otimes y)) = g(M_{AB^*}u_1,u_2).
\end{equation}

\begin{proof}
$(T_G(u_1\otimes x), T_G(u_2\otimes y))=(T_G^2(u_1\otimes x), u_2\otimes y)=
g(T_G(u_1\otimes x), u_2\otimes y)$  
(for properties of $T_G$, see remark~\ref{rem:pot}).
\[
\left(T_G(u_1\otimes x), u_2\otimes y\right)=
\left(\sum_{\sigma\in G}M(\sigma)u_1\otimes P(\sigma)x, u_2\otimes y\right) = 
\]
\[
\sum_{\sigma\in G}(M(\sigma)u_1, u_2)(P(\sigma)x, y) =
\sum_{\sigma\in G}(M(\sigma)u_1, u_2) \prod_{i=1}^n (x_{\sigma^{-1}(i)},y_i)=
\]
\[
\sum_{\sigma\in G}(M(\sigma)u_1, u_2) \prod_{i=1}^n AB^*(\sigma^{-1}(i), i)\,\,({\rm equation}\,\,
\ref{eq:ipmtrx})
\]
This latter equation becomes
\[
\left(\left(\sum_{\sigma\in G}M(\sigma)\prod_{i=1}^n AB^*(\sigma^{-1}(i), i)\right)u_1, u_2\right) =
\left(M_{AB^*}u_1, u_2\right).
\]
\end{proof}
\end{lemma}
%END LEM
%THM
\begin{theorem}[\bfseries Marcus's generalization of Schur's inequality \cite{mm:crs}]
\label{thm:mgsi}
Let $A$ and $B$ be $n\times n$ complex matrices.  
Let $M_K$ be the generalized matrix function of $K$. 
%as in definition~\ref{def:sgmf}.
Then, for any $u_1, u_2\in U$, 
\begin{equation}
\label{eq:mgsi}
\left| \left(M_{AB^*}u_1,u_2\right)\right|^2 \leq 
 \left(M_{AA^*}u_1,u_1\right) \left(M_{BB^*}u_2,u_2\right).
\end{equation}
\begin{proof}
Let $T_G(u_1\otimes x)$ and $T_G(u_2\otimes y)$ be generalized symmetric tensors where 
$x= x_1\otimes \cdots \otimes x_n$ and $y=y_1\otimes \cdots \otimes y_n$
are chosen so the the $n\times n$ matrices
$A$ and $B$ are associated with  $x$ and $y$ respectively
(remark~\ref{rem:ippm}).  
By the Cauchy-Schwarz inequality, we have
\begin{equation}
\label{eq:csi}
|\left(T_G(u_1\otimes x), T_G(u_2\otimes y)\right)|^2\leq
\end{equation}
\[
\left(T_G(u_1\otimes x), T_G(u_1\otimes  x)\right)\; 
\left(T_G(u_2\otimes y), T_G(u_2\otimes y)\right).
\]
From lemma~\ref{lem:ipgst}, equation~\ref{eq:csi} becomes
\begin{equation}
\label{eq:csimf}
g^2|(M_{AB^*}u_1, u_2)|^2 \leq g\,(M_{AA^*}u_1, u_1)\,\,g\,(M_{BB^*}u_2, u_2).\\
\end{equation}
\end{proof}
\end{theorem}
%END THM

%COR
\begin{corollary}[\bfseries Schur's inequality from theorem~\ref{thm:mgsi}]
\label{cor:mpsi}
Let $G$ be a subgroup of the symmetric group $S_n$ of degree $n$. 
Let $H=(h_{ij})$ be an $n\times n$ positive definite Hermitian matrix. 
Let $M$ be a representation of $G$ as unitary linear operators on $U$,
$\dim(U)=m$, and   
let $M_H = \sum_{\sigma\in G} M(\sigma) \prod_{i=1}^n h_{i\sigma(i)}$.
Then $M_H$ is a positive definite Hermitian transformation on $U$ and, if $u\in U$ has  $||u||=1$, then
\begin{equation}
\label{eq:sinc}
\det(H)\leq (M_Hu,u).
\end{equation}
\begin{proof}
Write 
$M_H=\sum_{\sigma\in G} M(\sigma) \prod_{i=1}^n H(i,\sigma(i))$.
\[
(M_H)^* = \sum_{\sigma\in G} (M(\sigma))^* \prod_{i=1}^n \overline{H(i,\sigma(i))}
= \sum_{\sigma\in G} (M(\sigma^{-1})) \prod_{i=1}^n H^*(\sigma(i),i) = M_{H^*}.
\]

Thus, $H^*=H$ implies that $(M_H)^* = M_H$  and hence $M_H$ is Hermitian.
If $H$ is positive definite Hermitian then $\det(H)>0$, and \ref{eq:sinc} will imply 
that every eigenvalue of $M_H$ is positive, hence $M_H$ is positive definite Hermitian. To prove~\ref{eq:sinc}, take $H=AA^*$ and let $B=I_n$ in equation~\ref{eq:mgsi}.  
 In which case $M_{BB^*}=M_{I_n} = I_m$.  
 Assume without loss of generality that $A$ is triangular and take 
 $u= u_1=u_2$ to be a unit vector.  In this case, $(M_Au,u) = \det(A)$ and
 equation~\ref{eq:mgsi} becomes $|\det(A)|^2 = \det(H) \leq (M_H u, u)$
 which proves~\ref{eq:sinc}. 
\end{proof}
\end{corollary}
%END COR

%SECTION
\section{Equalities}

The proof of Schur's inequality, corollary~\ref{cor:mpsi}, was obtained from Marcus's tensor form of the Cauchy-Schwartz inequality
\begin{equation}
\label{eq:csi2}
|\left(T_G(u_1\otimes x), T_G(u_2\otimes y)\right)|^2\leq
\end{equation}
\[
\left(T_G(u_1\otimes x), T_G(u_1\otimes  x)\right)\; 
\left(T_G(u_2\otimes y), T_G(u_2\otimes y)\right)
\]
by taking $u_1=u_2=u$ and 
$y=y_1\otimes  \cdots \otimes y_n = e_1\otimes  \cdots \otimes e_n = e$:
\begin{equation}
\label{eq:csi3}
|\left(T_G(u\otimes x), T_G(u\otimes e)\right)|^2\leq\|T_G(u\otimes x)\|^2\,|T_G(u\otimes e)\|^2.
\end{equation}
Both $T_G(u\otimes x)$ and $T_G(u\otimes e)$ are nonzero, hence equality occurs in equation~\ref{eq:csi3} if and only if 
\begin{equation}
\label{eq:csi4}
T_G(u\otimes x) = k T_G(u\otimes e)\,, \;k\neq 0.\\
\end{equation}\\
%REM
 \begin{remark}[\bfseries The case $T_G(u\otimes x) = k T_G(u\otimes e)\,, \;k\neq 0$, $n=2$ ]
 \label{rem:case2}
 Choose $u\in U,\, \|u\|=1.$ 
 The nonsingular upper triangular matrix  $A=(a_{ij})=((x_i,e_j))$ is in $\mathbf{M}_2(\bC)$.
 In this case, $x_1=a_{11}e_1 + a_{12}e_2$, $x_2=a_{22}e_2$ and  
 $u\otimes x = u\otimes x_1\otimes x_2$.  
 From \ref{def:gsop}
 \[
T_G(u\otimes x) = 
\sum_{\sigma\in G} M(\sigma)u\otimes  x_{\sigma^{-1}(1)} \otimes x_{\sigma^{-1}(2)}.
 \]
We have assumed that
\begin{equation}
\label{eq:parcastwo1}
T_G(u\otimes x) = k T_G(u\otimes e)\,, \;k\neq 0.
\end{equation}   
If $G=\{\iota, \tau\}$, $\tau = (1,2)$, we compute directly that
 \begin{equation}
 \label{eq:simpcasex}
 T_G(u\otimes x) = 
 a_{11}a_{22}\left(u\otimes e_1\otimes e_2 + M(\tau)u\otimes e_2\otimes e_1\right) + 
  a_{12}a_{22}Su\otimes e_2\otimes e_2
 \end{equation}
 where 
 $ 
S=I_m+M(\tau).
 $
 For $x=e$ (for which $A=I_2$) we get
 \begin{equation}
 \label{eq:simpcasee}
 T_G(u\otimes e) = u\otimes e_1\otimes e_2 + M(\tau)u\otimes e_2\otimes e_1.
 \end{equation}
 Let $\bG_A$ be as in remark~\ref{rem:comequcon}(3).
 From equations~\ref{eq:simpcasex} and  \ref{eq:simpcasee} we see that 
 equation~\ref{eq:parcastwo1} holds (for $G=\{\iota, \tau\}$) if and only if either 
\[
(1)\,\bG_A=\{\iota\}\;\;(a_{12}=0)\;\; {\rm or}  \;\;(2)\, \bG_A=\{\iota, \tau\} \;\;(a_{12}\neq 0)\;\;{\rm and}\;\;Su= 0. 
\]
SInce $S=I_m+M(\tau)$,  we have $(Su,u)=0$ if and only if $(M(\tau)u, u) = -1$.
 In the trivial case where $G=\{\iota\}$, equation~\ref{eq:parcastwo1} holds if and only if
 $\bG_A=\{\iota\}$ (i.e., $a_{12}=0$): 
 \begin{equation}
 \label{eq:simpcasex2}
 T_G(u\otimes x) = 
 a_{11}a_{22}u\otimes e_1\otimes e_2 +  a_{12}a_{22}u\otimes e_2\otimes e_2
 \end{equation}
 and
  \begin{equation}
 \label{eq:simpcasee2}
 T_G(u\otimes e) = u\otimes e_1\otimes e_2.
 \end{equation}
 To summarize, in the $n=2$ case, equation~\ref{eq:parcastwo1} holds if and only if
 $\bG_A$ is contained in $G$ and $(M(\sigma)u,u)=\epsilon(\sigma)$ for $\sigma\in \bG_A$.\\
  \end{remark}
\begin{definition}[\bfseries Compatible permutations]
\label{def:comper}
Let $\alpha\in \Gamma_n$.  A permutation $\sigma\in S_n$ will be
$\alpha$-compatible if $\alpha(\sigma(i))\geq i$ for all $i\in \underline{n}$.
Let
\[
S_n^\alpha = \{\sigma\mid \sigma\in S_n, \alpha(\sigma(i))\geq i, i\in\underline{n}\}
\]
denote the set of all $\alpha$-compatible  permutations.\\
\end{definition}
%DEF
%\begin{definition}[\bfseries Spike functions]
%\label{def:spifun}
%Let $r,c\in{\underline{n},\;r<c}$.  
%The $r,c$~{\em spike function}, $\alpha_{rc}$, is defined by
%\[
%\alpha_{rc}(r) =c\;{\rm and}\; \alpha_{rc}(i) =i,\; i\neq r.
%\]\\
%\end{definition}
%DEF
\begin{definition}[\bfseries Restricted $\alpha$-compatible permutations]
\label{def:resalpcomper}
$T_G(u\otimes x)$ has
upper triangular matrix $A=(a_{ij})=((x_i,e_j))$ associated with $x=x_1\otimes \cdots \otimes x_n$
and the basis $\{e_i\mid i=1, \ldots n\}.$ 
Define
\begin{equation}
\label{eq:resalpcomper}
S_n^\alpha (A)\equiv S_n^\alpha (x,e)=
\{\sigma \mid \sigma\in S_n^\alpha,\; \prod_{i=1}^n a_{i\alpha\sigma(i)}\neq 0\}.
\end{equation}
We call $S_n^\alpha (A)\equiv S_n^\alpha (x,e)$ the $\alpha$-compatible permutations restricted by $A$
or, alternatively,  
the $\alpha$-compatible permutations restricted by 
$x=x_1\otimes\cdots\otimes x_n$ and the orthonormal basis 
$e_i$, $1\leq i\leq  n$.\\
\end{definition}

%REM
\begin{remark}[\bfseries Inner products and $\alpha$-compatible permutations]
\label{rem:innprocomper}
Let $u\in U,\, \|u\|=1.$
Let $T_G(u\otimes x)$ be as in definition~\ref{def:resalpcomper}. 
Then
\begin{equation}
\left(T_G(u\otimes x), u\otimes e_\alpha\right)=
\sum_{\sigma\in G} \left(M(\sigma)u,u\right)\prod_{i=1}^n \left(x_i, e_{\alpha\sigma(i)}\right).
\end{equation}
Note that $A$ upper triangular implies that 
$\prod_{i=1}^n \left(x_i, e_{\alpha\sigma(i)}\right) = 0$ 
if $\sigma\not\in S_n^\alpha$ (definition~\ref{def:comper}).
In fact, $\prod_{i=1}^n \left(x_i, e_{\alpha\sigma(i)}\right) = 0$ if 
$\sigma\not\in S_n^\alpha(x,e)\subseteq S_n^\alpha $
(definition~\ref{def:resalpcomper}) and we have
\begin{equation}
\label{eq:spiforinnpro}
\left(T_G(u\otimes x), u\otimes e_\alpha\right)=
\sum_{\sigma\in G\cap S_n^\alpha (x,e)} \left(M(\sigma)u,u\right)
\prod_{i=1}^n \left(x_i, e_{\alpha\sigma(i)}\right).\\
\end{equation}
\end{remark}
%END REM
 
 We now discuss the general case of equality as it relates to collinearity:
\begin{equation}
\label{eq:gencascol}
T_G(u\otimes x) = k T_G(u\otimes e)\,, \;k\neq 0
\end{equation}
The proofs will follow the ideas and notation developed in remark~\ref{rem:case2}. 

%REM
\begin{remark}[\bfseries Terminology and notation]
\label{rem:ternot}  
If $K\subseteq \underline{n}$ then $S_K={\rm PER}(K)$ denotes all permutations of $K$.
Thus, $S_n=S_{\underline{n}}$. 
The set of all functions from $K$ to $K$ is denoted by  
$\Gamma_K$.
Thus, $\Gamma_n = \Gamma_{\underline{n}}$.
For $p\in \underline{n}$. 
\[
\Gamma_{n,p} = \{\gamma\mid \gamma\in \Gamma_n, \gamma(p)=p, 
\gamma(i)\neq p\; {\rm if}\; i\neq p\}.
\]
Let $G_p=\{\sigma\mid \sigma(p)=p\}$ denote the stability subgroup of $G$ at $p$.
Let $K_p = \underline{n}\setminus \{p\} = \{1, \ldots, p-1\}\cup \{p+1, \ldots n\}$.
Let $G'_p = G_p|K_p$ be the restriction of the stability subgroup $G_p$ to $K_p$.
For $\gamma\in \Gamma_{n,p}$ let $\gamma'_p = \gamma | K_p$ and note that
$\Gamma_{K_p}=\{\gamma'_p\mid \gamma \in \Gamma_{n,p}\}$.
Also, $|\Gamma_{n,p}|=|\Gamma_{K_p}|={(n-1)}^{n-1}$, the map $\gamma\mapsto \gamma'_p$
providing the canonical bijection.

If the tensor $x=x_1\otimes\ldots\otimes x_n$ has associated matrix $A=(a_{ij})=((x_i, e_j))$, then  
define the tensor $x^p$ and vectors $x^p_i$, $i\in K_p$,
to have associated  matrix $A(p\mid p)$ in the same sense. 
Here we use the standard notation which defines $A(p\mid p)$ to be the submatrix of $A$ gotten by deleting  row $p$ and column  $p$ from $A$.
Thus,
\[
x^p=x^p_1\otimes \cdots\otimes x^p_{p-1}\otimes x^p_{p+1} \otimes\cdots \otimes x^p_n.
\]
Let $T_G=\sum_{\sigma\in G} M(\sigma)\otimes P(\sigma)$ be a generalized 
symmetry operator on $U\otimes [\otimes^n V]$ where $\dim(V)=n$ and $\dim(U)=m$.  
Note that $T_G$ has {\em degree} $m$ and {\em order} $n$ 
(definition~\ref{def:gsop}).
$V$ has orthonormal basis $\{e_i\mid 1\leq i\leq n\}$ and $\otimes^n V$ has
orthonormal basis $\{e_\alpha \mid \alpha \in \Gamma_n\}$ where 
$e_\alpha = e_{\alpha(1)} \otimes \cdots \otimes e_{\alpha(n)}$.
Let $T_{G'_p} = \sum_{\sigma\in G'_p} M(\sigma)\otimes P(\sigma)$ act on $U\otimes [\otimes^{n-1} V'_p]$.
$V'_p$ has orthonormal basis $\{e_i\mid i\in K_p\}$ and
$\otimes^{n-1} V'_p$  has orthonormal basis 
$\{e_{\gamma'_p} \mid {\gamma'_p} \in \Gamma_{K_p}\}$ where
$e_{\gamma'_p} = \otimes_{i\in K_p}e_{\gamma'_p(i)}$.
$T_{G'_p}$ has degree $m$ and order $n-1$.\\
\end{remark}
%END
%REM
\begin{remark}[\bfseries Maximal row spike functions]
\label{rem:jump}
In the following discussion we will use the  notions of 
{\em spike functions} (definition~\ref{def:spifun}),
 {\em maximal row spike  functions},  $\alpha_{rc}$ (definition~\ref{def:maxrowspifun}), 
 and the characterization of permutations compatible with maximal row spike functions
 (lemma~\ref{lem:rescomperspifun}).  
The statements of these two definitions and the lemma are all that will be needed to understand the proofs that follow.  
 In fact, lemma~\ref{lem:rescomperspifun} can easily be proved directly from these two definitions.  This material has been included in a separate section in order to state and prove lemma~\ref{lem:spifuncomper} which is the natural combinatorial setting for this discussion.
\end{remark}
\begin{lemma}[\bfseries If columns $A^{(c)}$, $c>1$, have $r<c$ with $a_{rc}\neq 0$]
\label{lem:allcolspifnc}
Choose $u\in U$, $\|u\|=1$. Suppose, for some $k\neq 0$, $T_G(u\otimes x) = k T_G(u\otimes e)$, and, for all $1<c\leq n$, 
column $A^{(c)}$ has some $r<c$ with $a_{rc}\neq 0$. Then
$\bG_A=G=S_n$ and $(M(\sigma)u,u)=\epsilon(\sigma)$ for all $\sigma$.
%PRF  
\begin{proof}
The condition on the columns $A^{(c)}$ implies that every column ($c>1$) has a {\em maximum row spike function}
 $\alpha_{rc}$ (definition~\ref{def:maxrowspifun}). 
Consider the equation
\begin{equation}
\label{eq:spiforide}
\left(T_G(u\otimes x), u\otimes e_{\alpha_{rc}}\right) = 
k \left(T_G(u\otimes e), u\otimes e_{\alpha_{rc}}\right).
\end{equation}
By equation~\ref{eq:spiforinnpro}, 
\begin{equation}
\label{eq:spiforinnpro2}
 \left(T_G(u\otimes x), u\otimes e_{\alpha_{rc}}\right)=
\sum_{\sigma\in G\cap S_n^{\alpha_{rc}} (x,e)} \left(M(\sigma)u,u\right)\prod_{i=1}^n \left(x_i, e_{\alpha_{rc}\sigma(i)}\right).
\end{equation}
By lemma~\ref{lem:rescomperspifun}, $S_n^{\alpha_{rc}} (x,e) = \{\iota, \tau\}$ where
$\iota$ is the identity permutation and $\tau=(r,c)$ is a transposition.
Thus, $\left(T_G(u\otimes x), u\otimes e_{\alpha_{rc}}\right)=$
\begin{equation}
\label{eq:spiforinnpro3}
(I_mu,u)\left(\prod_{i\neq r} (x_i,e_i)\right)(x_r,e_c)\;\;{\rm if}\;\;\tau\not\in G\;\;{\rm and}
\end{equation}
\begin{equation}
\label{eq:spiforinnpro4}
\left[(I_mu,u)+(M(\tau)u,u)\right]\left(\prod_{i\neq r} (x_i,e_i)\right)(x_r,e_c)\;\;{\rm if}\;\;\tau\in G.
\end{equation}
Note that by the definition of the maximal row spike function,  $\alpha_{rc}$, we have
$(x_r,e_c)=a_{rc}\neq 0$.
However, $\left(T_G(u\otimes e), u\otimes e_{\alpha_{rc}}\right)=0$.
Thus, equation~\ref{eq:spiforinnpro4}, not equation~\ref{eq:spiforinnpro3}, must hold, and
$\left[(u,u)+(M(\tau)u,u)\right] = 0$ or $(M(\tau)u,u)=\epsilon(\sigma)$.
By hypothesis, for every column, $A^{(c)}$, $c>1$, there is a maximum row spike  function 
 $\alpha_{rc}$.  
Thus, we have shown that for every $c$, $1<c\leq n$, there is an $r<c$ such that 
$\tau=(r,c)\in \bG_A$ satisfies $\left[(u,u)+(M(\tau)u,u)\right] = 0$ or $(M(\tau)u,u)=\epsilon(\sigma)$.  
This set of transpositions, $\tau$, generates $S_n$ (a trivial induction).  
We have proved that $(M(\tau)u,u)=\epsilon(\sigma)$ for all $\sigma\in \bG_A=S_n=G$
(note remark~\ref{rem:comequcon} here).\\
\end{proof}  
\end{lemma}
%END LEM

\begin{lemma}[\bfseries General case where column $A^{(p)}$ has only $a_{pp}\neq0$]
\label{lem:procasgen}
Choose $u\in U$, $\|u\|=1$. Assume for $k\neq 0$, $T_G(u\otimes x) = k T_G(u\otimes e)$.
Using the terminology of remark~\ref{rem:ternot}, assume that $x$ is associated with
 the upper triangular and nonsingular matrix $A=(a_{ij})$.  
Assume for some $p$, $1<p\leq n$, the only nonzero entry in column $A^{(p)}$ is the diagonal entry $a_{pp}$.
Then, for any $\gamma\in \Gamma_{n,p}$, 
\begin{equation}
 \label{eq:procasgen}
 \left(T_G(u\otimes x), u\otimes e_\gamma\right)=
 a_{pp}\left(T_{G'_p}(u\otimes x^p), u\otimes e_{\gamma'_p}\right).
 \end{equation}
 %PRF
 \begin{proof}
 Consider $u\otimes e_\gamma$, $\gamma\in \Gamma_{n,p}$.
Then,
\begin{equation}
 \label{eq:procasn1}
 \left(T_G(u\otimes x), u\otimes e_\gamma\right)=
 \sum_{\sigma\in G}(M(\sigma)u,u)\prod_{i=1}^n(x_{\sigma^{-1}(i)}, e_{\gamma(i)})
 \end{equation}
 where 
 \[
 \prod_{i=1}^n(x_{\sigma^{-1}(i)}, e_{\gamma(i)}) = 
 (x_{\sigma^{-1}(1)}, e_{\gamma(1)})\cdots (x_{\sigma^{-1}(p)}, e_{\gamma(p)})\cdots (x_{\sigma^{-1}(n)}, e_{\gamma(n)}).
 \]
 Note that $(x_{\sigma^{-1}(p)}, e_{\gamma(p)})=(x_{\sigma^{-1}(p)}, e_p)\neq 0$
 requires $\sigma^{-1}(p)=p$.   
 Thus, $(x_{\sigma^{-1}(p)}, e_p)=(x_p, e_p)=a_{pp}\neq 0.$
 Hence the sum in equation~\ref{eq:procasn1} can be taken over $G_p$,  
 the  stabilizer subgroup of $G$ at $p$.
 We have 
 \begin{equation}
  \label{eq:procasn3}
 \left(T_G(u\otimes x), u\otimes e_\gamma\right)=
 a_{pp}\sum_{\sigma\in G_p}(M(\sigma)u,u)\prod_{i\neq p} (x_{\sigma^{-1}(i)}, e_{\gamma(i)}).
 \end{equation}
 
As in remark~\ref{rem:ternot}, let $G'_p$ be $G_p$ restricted to $K_p$, and
 let $x^p=x^p_1\otimes \ldots x^p_{p-1}\otimes x^p_{p+1} \otimes x^p_n$
be the tensor that has matrix $A(p\mid p)$ with respect to the orthonormal basis
$\{e_i\mid i\in K_p\}$.
Let  $\gamma'_p$ denote $\gamma$ restricted to $K_p$.  

 The unitary representation $M$ of $G$ restricts in the obvious way to $G_p$ and ${G'}_p$.
 Let $V'_p = \langle e_i\mid i\in K_p\rangle$.
 The generalized symmetry operator 
 $T_{G'_p} = \sum_{\sigma\in G'_p}M(\sigma)\otimes P(\sigma)$ 
 acting on  $U\otimes [\otimes V'_p]$ allows us to reformulate equation~\ref{eq:procasn3}
 as follows:
 \begin{equation}
 \label{eq:procasgen2}
 \left(T_G(u\otimes x), u\otimes e_\gamma\right)=
 a_{pp}\left(T_{G'_p}(u\otimes x^p), u\otimes e_{\gamma'_p}\right).\\ 
\end{equation}
\end{proof}
\end{lemma}

\begin{remark}[\bfseries Inductive step using equation~\ref{eq:procasgen}]
\label{rem:sumprocasgen}
From  the terminology of remark~\ref{rem:ternot} we note that 
$T_{G'_p} = \sum_{\sigma\in G'_p} M(\sigma)\otimes P(\sigma)$ acts on $U\otimes [\otimes^{n-1} V'_p]$.
$V'_p$ has orthonormal basis $\{e_i\mid i\in K_p\}$ and
$\otimes^{n-1} V'_p$  has orthonormal basis 
$\{e_{\gamma'_p} \mid {\gamma'_p} \in \Gamma_{K_p}\}$ where
$e_{\gamma'_p} = \otimes_{i\in K_p}e_{\gamma'_p(i)}$.
Thus, $T_{G'_p}$ has degree $m$ and order $n-1$ and will be used in the inductive step (on $n$) in the next theorem.
\end{remark}

\begin{theorem}
\label{thm:mainresult}
Let $u\in U$, $\|u\|=1$. 
Let $A=(a_{ij})=((x_i,e_j))$ 
be the upper triangular nonsingular matrix associated with 
$x=x_1\otimes \cdots \otimes x_n$  and basis $e_i$, $1\leq i\leq n$.
If $\,T_G(u\otimes x) = k T_G(u\otimes e)\,, \;k\neq 0$, then  $\bG_A$ is a 
subgroup of $G$ and $(M(\sigma)u,u)=\epsilon(\sigma)$ for all $\sigma\in \bG_A$.
\begin{proof}
If for all $c$, $1<c\leq n$, column $A^{(c)}$ has some $r<c$ with $a_{rc}\neq 0$
then the result follows from lemma~\ref{lem:allcolspifnc}. 
Otherwise, for some $p$, $1<p\leq n$, the only nonzero entry in column $A^{(p)}$ is the diagonal entry $a_{pp}$.
We apply lemma~\ref{lem:procasgen}, in particular
 \begin{equation}
 \label{eq:mainresult1}
 \left(T_G(u\otimes x), u\otimes e_\gamma\right)=
 a_{pp}\left(T_{G'_p}(u\otimes x^p), u\otimes e_{\gamma'_p}\right)
\end{equation}
and for $x=e$:
\begin{equation}
 \label{eq:mainresult2}
 \left(T_G(u\otimes e), u\otimes e_\gamma\right)=
 \left(T_{G'_p}(u\otimes e^p), u\otimes e_{\gamma'_p}\right).\\ 
\end{equation}
The proof is by induction on the order $n$ of $T_G$ (i.e., $n=\dim(V)$).
The base case, $n=2$, is established in  remark~\ref{rem:case2}.
Assume the theorem is true for generalized symmetry operators of degree $m$ and order $n-1$.
Note that $T_G(u\otimes x) = k T_G(u\otimes e)\,, \;k\neq 0$, implies that 
$(T_G(u\otimes x), u\otimes e_\gamma)=k(T_G(u\otimes e), u\otimes e_\gamma)$
for all $\gamma\in \Gamma_{n}$.  Thus, by equations~\ref{eq:mainresult1} and \ref{eq:mainresult2}
\begin{equation}
\label{eq:mainresult3}
a_{pp}(T_{G'_p}(u\otimes x^p), u\otimes e^p_{\gamma'_p})=
k(T_{G'_p}(u\otimes e^p), u\otimes e^p_{\gamma'_p})
\end{equation}
for all $\gamma'_p\in\Gamma_{K_p}$.  
But, $\{u\otimes e_{\gamma'_p} \mid {\gamma'_p} \in \Gamma_{K_p}\}$
spans $U\otimes [\otimes^{n-1}V_p]$ (\ref{rem:ternot}).  
Thus, 
\begin{equation}
\label{eq:mainresult4}
(T_{G'_p}(u\otimes x^p))= k_p(T_{G'_p}(u\otimes e^p)),\;\;k_p=k/a_{pp}.
\end{equation}
$T_{G'_p}$ has degree $m$ and order $n-1$.
Likewise,
\begin{equation}
\label{eq:mainresult5}
T_{G'_1}(u\otimes x^1) = k_1 T_{G'_1}(u\otimes e^1),\;\;k_1=k/a_{11}.
\end{equation}
$T_{G'_1}$ has degree $m$ and order $n-1$.
The tensor $x^1$ is associated with the 
upper triangular nonsingular submatrix $B=A(1|1)$.  
By \ref{eq:mainresult5} and the induction hypothesis, $\bG_B$ is a subgroup of $G'_1$, the stabilizer subgroup of $G$ at $1$ restricted to $K_1$, and $(M(\sigma)u,u)=\epsilon(\sigma)$ for $\sigma\in G_B$.
Let $C=A(p|p)$. $C$ is also upper triangular and nonsingular. 
By \ref{eq:mainresult4} and the induction hypothesis, $\bG_C$ is a subgroup of $G'_p$, the stabilizer subgroup of $G$ at $p$ restricted to $K_p$, and $(M(\sigma)u,u)=\epsilon(\sigma)$ for $\sigma\in \bG_C$.
Observe that  the union of the set of generating transpositions 
$\{(i,j)\mid i<j,\, B(i,j)\neq 0\}$, 
together with the set $\{(i,j)\mid i<j,\, C(i,j)\neq 0\}$, equals the set
$\{(i,j)\mid i<j,\, A(i,j)\neq 0\}$.  Here we use the fact that $a_{1p}= 0$.
Thus, $\bG_A$ is a subgroup of $G$ and $(M(\sigma)u,u)=\epsilon(\sigma)$ for all $\sigma\in G_A$.
We use remark~\ref{rem:comequcon}(2) here.
\end{proof}
\end{theorem}
%END LEM

%SECTION
\section{Remarks about compatible permutations}

We combine here some interesting combinatorial lemmas about compatible permutations.
We have used only  lemma~\ref{lem:rescomperspifun} which can be proved independently as an exercise.

\begin{definition}[\bfseries Spike functions]
\label{def:spifun}
Let $r,c\in{\underline{n},\;r<c}$.  
The $r,c$~{\em spike function}, $\alpha_{rc}$, is defined by
\[
\alpha_{rc}(r) =c\;{\rm and}\; \alpha_{rc}(i) =i,\; i\neq r.
\]\\
\end{definition}
%LEM
The following lemma characterizes the spike function compatible permutations 
(\cite{gw:tch}, Lemma 3.3, p. 339).\\

\begin{lemma}[\bfseries Spike function compatible permutations]
\label{lem:spifuncomper}
Let $\alpha_{rc}$ be a spike function (\ref{def:spifun}).
The set of all $\alpha_{rc}$-compatable permutations, $S_n^{\alpha_{rc}}$, consists of the identity permutation, $\iota$, together with all $\sigma$ which have only one cycle of length greater than one.   
Moreover, that cycle can be written in the form  $(r_1,\ldots , r_p)$ where 
$r=r_1 < \cdots < r_p\leq c$. 
Thus, $|S_n^{\alpha_{rc}}|=2^{c-r}$.
%PROOF
\begin{proof}
Let $(r_1,\ldots , r_p)$ be a nontrivial cycle of $\sigma$. 
Assume without loss of generality that $r_1 = \min\{r_1,\ldots , r_p\}$.
Suppose $r_1\in \{1, \ldots, r-1, r+1, \ldots n\}$ (i.e., $r_1\neq r$).
Then $\alpha_{rc}(\sigma(r_p)) = \alpha_{rc}(r_1) = r_1 < r_p$ which
contradicts the assumption $\sigma\in S_n^{\alpha_{rc}}$.
Thus, $r_1=r$ which implies there is only one nontrivial cycle of $\sigma$. 
If $r_p>c$, then 
$\alpha_{rc}(\sigma(r_p)) = \alpha_{rc}(r) = c < r_p$, again contradicting
$\sigma\in S_n^{\alpha_{rc}}$.
Thus, $r_p\leq c$ and we need only show the strictly increasing property:
$r=r_1 < \cdots < r_p$.  By minimality of $r_1$, we have $r_1< r_2$.   
Suppose we have $r_1< r_2 < \cdots < r_{t-1} > r_t$ for some $2< t \leq p$.
Then $\alpha_{rc}(\sigma(r_{t-1})) = \alpha_{rc}(r_t) = r_t < r_{t-1}$
which contradicts the assumption $\sigma\in S_n^{\alpha_{rc}}$.
Thus the unique nontrivial cycle of $\sigma$ is of the form
$(r_1,\ldots , r_p)$ which we may assume satisfies $r=r_1 < \cdots < r_p\leq c$. 
These cycles can be constructed by choosing a nonempty subset of
$\{r+1, \ldots, c\}$.  This can be done in $2^{c-r} - 1$ ways.  Counting the 
identity permutation, this gives $|S_n^{\alpha_{rc}}|=2^{c-r}$.\\
\end{proof}
\end{lemma}
%END LEMMA
\begin{definition}[\bfseries Maximal row spike functions]
\label{def:maxrowspifun}
We use the terminology of remark~\ref{rem:innprocomper}.
Let $A=(a_{ij}) = ((x_i, e_j))$ be the upper triangular nonsingular matrix associated with  $x_1\otimes\dots\otimes x_n$ and basis $e_1, \ldots, e_n$. 
Suppose column $A^{(c)}$ has at least two nonzero entries.
Let $r=\max\{i\mid i<c, a_{ic} \neq 0\}$.  
We call the spike  function $\alpha_{rc}$ the {\em maximal row} spike function for column $c$.  If  $\{i\mid i<c, a_{ic} \neq 0\}$ is empty, we say column $c$ has no maximal row spike function.
\end{definition}

We next characterize the sets $S_n^\alpha(x,e)$ for maximal row spike functions 
$\alpha = \alpha_{rc}$.\\
%LEM
\begin{lemma}[\bfseries Permutations compatible with maximal row spike functions]
\label{lem:rescomperspifun}
We use definition~\ref{def:maxrowspifun} and the terminology of 
remark~\ref{rem:innprocomper}.
Let $A=(a_{ij})=((x_i,e_j))$ be the upper triangular matrix associated with $x=x_1\otimes \cdots \otimes x_n$
and the orthonormal basis $\{e_i\mid i=1, \ldots n\}.$ 
Let $\alpha_{rc}$ be a maximal row spike function for column $A^{(c)}$.
Then
\begin{equation}
S_n^{\alpha_{rc}}(x,e) = \{\iota, \tau\}
\end{equation}
where $\iota$ is the identity and $\tau=(r,c)$.
%PRF
\begin{proof}
Note that $S_n^\alpha(x,e)\subseteq S_n^{\alpha_{rc}}$.  
By lemma~\ref{lem:spifuncomper}, $\sigma\in S_n^{\alpha_{rc}}$ implies
$\sigma=(r_1,\ldots , r_p)$ where  $r=r_1 < \cdots < r_p\leq c$. 
If $r_p<c$ then 
\[
A(r_p, \alpha_{rc}(\sigma(r_p)) = A(r_p, \alpha_{rc}(r)) = A(r_p, c) = 0
\]
 by maximality of $r$.
Thus, by equation~\ref{eq:resalpcomper}, $\sigma\not\in S_n^{\alpha_{rc}}(x,e).$
Thus, $r_p=c$.  
If $p>2$, we have $1< r_{p-1} < c=r_p$.  
\[
A(r_{p-1}, \alpha_{rc}(\sigma(r_{p-1})) = A(r_{p-1}, \alpha_{rc}(c)) = A(r_{p-1}, c) = 0
\]
by the maximality of $r$.
Again by equation~\ref{eq:resalpcomper}, $\sigma\not\in S_n^{\alpha_{rc}}(x,e).$
Thus, $p=2$.    
We have shown that  $S_n^{\alpha_{rc}}(x,e) = \{\iota, \tau\}$ where
$\iota$ is the identity and $\tau=(r,c)$.\\
\end{proof}
%END PRF
\end{lemma}
\section{Character form and examples of Schur's Theorem}
The following definition and terminology will be useful:
\begin{definition} [\bfseries Schur generalized matrix function]
\label{def:sgmf}
Let $G$ be a subgroup of the symmetric group $S_n$ of degree $n$. 
Let $K=(k_{ij})$ (alternatively, $k_{ij}\equiv K(i,j)$) be an $n\times n$ matrix with complex entries.
Let $M$ be a representation of $G$ as unitary linear operators on a unitary space $U$, 
$\dim U=m$.
Then the {\em Schur generalized matrix function}, $M_K$, is defined as follows:
\begin{equation}
\label{eq:stmf}
M_K = \sum_{\sigma\in G} M(\sigma) \prod_{i=1}^n K(\sigma^{-1}(i), i) =
\sum_{\sigma\in G} M(\sigma) \prod_{i=1}^n K(i,\sigma(i)).
\end{equation}
\end{definition}
\hspace*{1in}\\

%REMARK CHARACTER FORM
\begin{remark}[\bfseries Character or trace form of Schur's theorem]
\label{rem:cfst} 
We use the notation of theorem~\ref{thm:stf}.
$M_H$ is positive definite Hermitian.  
Let $\mathcal{E}_{M_H}$ denote the multiset of eigenvalues of $M_H$ (i.e., eigenvalues with multiplicities).
We have $\lambda >0$ for all $\lambda\in\mathcal{E}_{M_H}$ and
\begin{equation}
\label{eq:trcform}
\sum_{\lambda\in\mathcal{E}_{M_H}} \lambda = {\rm Tr}(M_H) = 
\sum_{\sigma\in G}{\rm Tr}(M(\sigma))\prod_{i=1}^n h_{i\sigma(i)}.
\end{equation}
Let $u_{\rm min}$ be a unit eigenvector corresponding to $\lambda_{\rm min}$,
the minimum eigenvalue in $\mathcal{E}_{M_H}$.
From theorem~\ref{thm:stf} we have 
$\det(H)\leq (M_Hu_{\rm min},u_{\rm min})=\lambda_{\rm min}$.
From equation~\ref{eq:trcform}, $m\det(H)\leq m\lambda_{\rm min}\leq {\rm Tr}(M_H)$.
Thus, $m\det(H) = {\rm Tr}(M_H)$ if and only if 
$\mathcal{E}_{M_H}=\{\det(H), \ldots, \det(H)\}$ (i.e., $\det(H)$ has multiplicity $m$).
In this case, $M_H=\det(H)I_m$, $I_m$ the identity.  
We have shown
\begin{equation}
\label{eq:trcform2}
m\det(H)\leq {\rm Tr}(M_H) \; {\rm and}\;
\end{equation}
\[
m\det(H) = {\rm Tr}(M_H) \iff M_H=\det(H)I_m.
\]
As an example, take $G=S_3$, the symmetric group on $\{1,2,3\}$, and take
\[
M(\sigma)=
\left(
\begin{array}{cc}
\epsilon(\sigma)&0\\
0&1
\end{array}
\right).
\]
Let $H$ be a $3\times 3$ positive definite Hermitian matrix such that the 
permanent, ${\rm per}(H) > \det(H)$.  Then
\[
M_H=
\left(
\begin{array}{cc}

\det(H)&0\\
0&{\rm per}(H)
\end{array}
\right).
\]
In this case,
\[
(M_Hu_{\rm min},u_{\rm min})=\lambda_{\rm min} = \det(H).
\] 
Note that the equality condition of \ref{eq:trcform2} does not hold:
$2 \det(H) < {\rm Tr}(M_H) = \det(H) + {\rm per}(H)$. 
To compare this with the equality condition~\ref{eq:seq}, note that  $\bG_H\subseteq S_3$ holds trivially.  
Let $u=(u_1,u_2)\in U$ have $\|u\|=1$.
Condition~\ref{eq:seq} states that $(M_Hu,u) = \det(H)$ if and only if
$(M(\sigma)u,u)=\epsilon(\sigma)$ for all $\sigma\in \bG_H$.
In this case these conditions hold if and only if $u=(e^{ir},0)$.
Note that ${\rm per}(H) > \det(H)$ implies $H$ has at least one off diagonal element.
%Note, $(M(\tau)u,u)=\epsilon(\tau)|u_1|^2 + |u_2|^2=\epsilon(\tau)=-1$  
%for any transposition  $\tau\in \bG_H$ implies that $u=(e^{ir},0)$ 
%(${\rm per}(H) > \det(H)$ implies $H$ has at least one off diagonal element). 
\end{remark}
%END REMARK

%BEGIN REMARK EXAMPLE SCHUR'S THEOREM
\begin{remark}[\bfseries Examples of Schur's theorem]
\label{rem:exscth}
We consider an example of $M_H = \sum_{\sigma\in G} M(\sigma) \prod_{i=1}^n h_{i\sigma(i)}$.
Take $G=S_3$, the symmetric group on $\{1,2,3\}$.  
Take the unitary representation to be the following:
\[
M(e)=\left(\begin{array}{cc}1&0\\0&1\end{array}\right)\;\;
M((123))=\left(\begin{array}{cc}\frac{-1}{2}&\frac{+\sqrt{3}}{2}\\\frac{-\sqrt{3}}{2}&\frac{-1}{2}\end{array}\right)\;\;
M((132))=\left(\begin{array}{cc}\frac{-1}{2}&\frac{-\sqrt{3}}{2}\\\frac{+\sqrt{3}}{2}&\frac{-1}{2}\end{array}\right)\;\;
\]
\[
M(23)=\left(\begin{array}{cc}\frac{-1}{2}&\frac{+\sqrt{3}}{2}\\\frac{+\sqrt{3}}{2}&\frac{1}{2}\end{array}\right)\;\;
M((13))=\left(\begin{array}{cc}\frac{-1}{2}&\frac{-\sqrt{3}}{2}\\\frac{-\sqrt{3}}{2}&\frac{1}{2}\end{array}\right)\;\;
M((12))=\left(\begin{array}{cc}1&0\\0&-1\end{array}\right)\;\;
\]

The Schur generalized matrix function, $M_H$ is
\[
M_H=
\left(\begin{array}{cc}1&0\\0&1\end{array}\right)
h_{11}h_{22}h_{33}\;\;+\;\;
\left(\begin{array}{cc}\frac{-1}{2}&\frac{+\sqrt{3}}{2}\\\frac{-\sqrt{3}}{2}&\frac{-1}{2}\end{array}\right)
h_{12}h_{23}h_{31}\;\;+\;\;
\left(\begin{array}{cc}\frac{-1}{2}&\frac{-\sqrt{3}}{2}\\\frac{+\sqrt{3}}{2}&\frac{-1}{2}\end{array}\right)
h_{13}h_{21}h_{32}\;\;+\;\;
\]
\[
\left(\begin{array}{cc}\frac{-1}{2}&\frac{+\sqrt{3}}{2}\\\frac{+\sqrt{3}}{2}&\frac{1}{2}\end{array}\right)
h_{11}h_{23}h_{32}\;\;+\;\;
\left(\begin{array}{cc}\frac{-1}{2}&\frac{-\sqrt{3}}{2}\\\frac{-\sqrt{3}}{2}&\frac{1}{2}\end{array}\right)
h_{13}h_{22}h_{31}\;\;+\;\;
\left(\begin{array}{cc}1&0\\0&-1\end{array}\right)
h_{12}h_{21}h_{33}.
\]
Let
\[
H=
\left(
\begin{array}{ccc}
1&0&0\\
0&3&i\\
0&-i&1
\end{array}
\right)
\;\;\;\;{\rm which\;\; gives\;\;}\;\;
M_H =
\begin{pmatrix}
%\left(
%\begin{array}{cc}
\frac{5}{2}&\frac{\sqrt{3}}{2}\\
\frac{\sqrt{3}}{2}&\frac{7}{2}
\end{pmatrix}
%\end{array}
%\right).
\] 
The eigenvalues $\mathcal{E}_{M_H}=\{2, 4\}$, $\det(H)=2<{\rm Trace}(M_H)=6$.
Also, 
\[
(M_Hu_{\rm min}, u_{\rm min}) = \lambda_{\rm min} = 2 = \det(H)
\]
By theorem~\ref{thm:stf}, equation~\ref{eq:seq},
$\det(H) = (M_Hu_{\rm min},u_{\rm min})$ implies that 
\[
\bG_H\subseteq G\,\, {\rm and}\,\,(M(\sigma)u_{\rm min},u_{\rm min})=
\epsilon(\sigma),\, \sigma\in \bG_H.
\]
By definition (in theorem~\ref{thm:stf}), $\bG_H$ is the group generated by 
$\{e, (2,3)\}$.
We have, trivially, $(M(e)u_{\rm min},u_{\rm min})=1$.   
To compute $(M(2,3)u_{\rm min},u_{\rm min})$ we evaluate the minimum eigenvector
$u_{\rm min}$ for $M_H$ and get $u_{\rm min} = (\frac{-\sqrt{3}}{2}, \frac{1}{2})$
and compute $M((23))u_{\rm min}$ where 
\[
M((23))=\left(\begin{array}{cc}\frac{-1}{2}&\frac{+\sqrt{3}}{2}\\\frac{+\sqrt{3}}{2}&\frac{1}{2}\end{array}\right)
\]
getting $-u_{\rm min}$ so $(M((23))u_{\rm min}, u_{\rm min})=-1=\epsilon((23))$ as required by the equality condition~\ref{eq:seq}.\\
\end{remark}

{\bf Acknowledgment:}  The author wishes to thank Dr. Tony Trojanowski for many helpful observations and specific suggestions that improved this paper.\\
%
%S. Gill Williamson, 2012\\
%${\rm cseweb.ucsd.edu\slash\sim gill\slash}$  
%
%\tableofcontents
%\index{contents}
%\hyperlink{index}{Index}
%\newpage
%\thispagestyle{empty}

\bibliographystyle{alpha}
\bibliography{SchurIneq}

%\label{cor:relpridiaent}
%\newpage
%\mbox{}
%\newpage
%\hypertarget{index}{ }
%\printindex
%\newpage
%\centerline{NOTES}
%\newpage
%\centerline{NOTES}
%\newpage
%\centerline{NOTES}
\end{document}